\documentclass[11pt]{amsart}
\usepackage{amsmath, amssymb, amsthm, mathtools, enumitem}
\usepackage[noabbrev,capitalize]{cleveref}
\usepackage[textheight=8.75in, textwidth=6.85in]{geometry}

\newtheorem{theorem}{Theorem}[section]
\newtheorem{lemma}[theorem]{Lemma}
\newtheorem{corollary}[theorem]{Corollary}
\newtheorem{proposition}[theorem]{Proposition}
\theoremstyle{remark}
\newtheorem{remark}{Remark}
\numberwithin{equation}{section}

\newcommand{\F}{\mathbb F}
\newcommand{\Z}{\mathbb Z}
\newcommand{\Q}{\mathbb Q}
\newcommand{\C}{\mathbb C}
\newcommand{\PP}{\mathbb P}
\newcommand{\Fpb}{\overline{\F}_p}
\newcommand{\Fpp}{\F_{p^2}}
\newcommand{\eps}{\varepsilon}
\newcommand{\lam}{\lambda}
\newcommand{\ord}{\operatorname{ord}}
\DeclareMathOperator{\Gal}{Gal}
\newcommand{\Aut}{\operatorname{Aut}}
\newcommand{\End}{\operatorname{End}}
\newcommand{\Frob}{\operatorname{Frob}}
\newcommand{\SL}{\mathrm{SL}}
\newcommand{\GL}{\mathrm{GL}}
\newcommand{\wt}[1]{\widetilde{#1}}

\title[Modular forms for chromatic homotopy: supersingular congruences]{Modular forms for chromatic homotopy: supersingular congruences}

\author{Ken Ono}
\dedicatory{In memory of Jack Morava}

\address{Dept. of Mathematics, University of Virginia, Charlottesville, VA 22904, USA}
\email{ko5wk@virginia.edu}

\begin{document}

\begin{abstract}
We prove a conjecture of Larson in Behrens' program on congruences of modular forms attached to the divided beta family in the Adams--Novikov spectral sequence for the stable homotopy groups of spheres. The conjecture gives a sharp criterion for when the modular form associated to a divided beta element can be represented by a pure power of the discriminant modular form. Writing $i=rp^{n}$ with $(r,p)=1$ and $t=i(p^2-1)/12$, Larson's conjecture asserts that the Behrens form $f_{i/j}$ (which is well defined modulo $p$) may be taken to be the pure power $\Delta^{t}$ precisely when $1\le j\le p^{n}$, and admits no such representative otherwise. We prove this for all primes $p\ge5$. The proof reduces the decisive congruence condition to a geometric statement on supersingular points of modular curves. Namely, that for every prime $\ell\ne p$, the value of the modular function $V_\ell(\Delta)/\Delta$ at each supersingular point of $X_0(\ell)$ is an $(p^2-1)/12$-th root of unity.
\end{abstract}

\maketitle

\section{Introduction and statement of results}\label{sec:intro}

Behrens' congruences attach modular forms modulo $p$ to order $p$ classes in the divided beta family in the Adams--Novikov spectral sequence (ANSS). We use only the congruence formulation of this construction, but recall the setting. At a fixed prime $p$, the ANSS is built from the cohomology of the moduli of $p$-adic formal groups and has the form
\[
\operatorname{Ext}^{\,s,t}_{BP_\ast BP}(BP_\ast,\,BP_\ast)\;\Longrightarrow\;\pi_{t-s}^{\,S},
\]
see \cite{Adams,Ravenel}. For odd $p$, its $s=2$ line contains the divided beta family, with elements denoted $\beta_{i/j}$ here. This family is one of the basic chromatic height $2$ families in the stable homotopy groups of spheres, and explicit descriptions of its representatives are a central testing ground for the relationship between stable homotopy theory and the geometry of elliptic curves.

Building on the elliptic cohomology viewpoint of Ando--Hopkins--Strickland \cite{AHS}, Behrens \cite{Behrens,BehrensCongruences} associated to each order $p$ divided beta element a modular form modulo $p$ satisfying four congruence conditions. In this form, the problem of identifying divided beta elements becomes a concrete problem about modular forms modulo $p$ and their behavior on modular curves. Behrens also suggested that, in many cases, these forms should be pure powers of
\[
\Delta(q):=\frac{E_4(q)^3-E_6(q)^2}{1728}=q\prod_{n=1}^{\infty}(1-q^n)^{24}.
\]
Here $q=e^{2\pi i\tau}$ with $\tau$ in the upper half-plane, and 
$$E_k=1-\tfrac{2k}{B_k}\sum_{n\ge1}\sigma_{k-1}(n)q^n
$$
is the normalized Eisenstein series of weight $k$, where $B_k$ is the $k$-th Bernoulli number and $\sigma_{k-1}(n)=\sum_{d\mid n}d^{k-1}$.

Larson \cite{Larson} made this prediction precise at level $2$. He rewrote the relevant divisibilities in terms of the Hasse invariant $A_p\equiv E_{p-1}\pmod p$ (cf.\ Katz \cite{Katz350}, Ono \cite{OnoCBMS}), proved the case $p=5$, exhibited examples for $p=7,11,13$, and $677$, and stated the conjecture proved below. Our main theorem completes this calculation for pure powers of $\Delta$ at every prime $p\ge5$: it gives the exact range in which these powers represent the Behrens forms attached to the divided beta family, and it proves that the boundary is sharp. The only topological input used in this paper is the existence of the Behrens forms $f_{i/j}$ and their four defining conditions (C1)--(C4), recalled in Section~\ref{sec:behrens}; after that point the proof is a calculation on modular curves and supersingular elliptic curves.

\medskip
\noindent
\textbf{Conjecture (Larson \cite[Conj.\ 5.21]{Larson}).}
\emph{For a prime $p\ge 5$, write $i=r p^{n}$ with $(r,p)=1$, and put $t=\dfrac{i(p^2-1)}{12}$. Then for every index $j$ admissible in the sense of Larson's Lemma~2.1 (see Section~\textup{\ref{sec:behrens}}),} we have that
\[
f_{i/j}\ =\ \Delta^{\,t}\qquad\text{for all admissible } j\text{ with } 1\le j\le p^n,
\]
\emph{and no representative of $f_{i/j}$ equals $\Delta^{\,t}$ for any other admissible value of $j$.
}

\smallskip
\begin{remark}
 Here $f_{i/j}$ is determined up to multiplication by a unit modulo $p$ and addition of a modular form  whose reduction modulo $p$ vanishes \textup{\cite[\S1]{Larson}}. Equivalently, on the range of admissible $j\le p^n$ the pure power $\Delta^{t}$ satisfies Behrens' conditions \textup{(C1)--(C4)}, while for admissible $j>p^{n}$ it fails them.
 \end{remark}
 
\medskip

We establish the conjecture in full.

\begin{theorem}\label{thm:main}
Let $p\ge 5$ be prime, write $i=rp^{n}$ with $(r,p)=1$, and set $t=\dfrac{i(p^2-1)}{12}$. Let $j$ be an index admissible for $\beta_{i/j}$ in the sense of Larson's Lemma~2.1 \textup{(see Section~\ref{sec:behrens})}.
\begin{enumerate}
\item[(i)] If $1\le j\le p^{n}$, then the pure power $\Delta^{t}$ satisfies Behrens' conditions \textup{(C1)--(C4)}. Namely, $f_{i/j}$ may be taken to be $\Delta^{t}$.
\item[(ii)] If $j>p^{n}$ (such admissible indices occur precisely when $r>1$ and $n\ge2$), then the pure power $\Delta^{t}$ fails condition \textup{(C4)}. Namely, no representative of $f_{i/j}$ equals $\Delta^{t}$.
\end{enumerate}
In particular, Larson's conjecture holds for all primes $p\ge5$.
\end{theorem}

\begin{remark}\label{rem:strongdiv}
The divisibility statements underlying Theorem~\ref{thm:main}, Propositions~\ref{prop:level2} and \ref{prop:bridge}  are established for \emph{every} integer $j\ge1$. Only the interpretation ``$\Delta^t$ represents $f_{i/j}$'' requires $(i,j)$ to be admissible, since only then is there a divided beta element $\beta_{i/j}$ and an attached Behrens form. The core of the proof is insensitive to Larson's exclusion of the indices $j=p,2p,\dots,b_{n-2}p$.
\end{remark}

Throughout the paper, we write
\begin{equation}
m:=\frac{p^2-1}{12}\in\Z,\qquad t=rp^nm,\qquad t':=\frac{t}{p^n}=rm,
\end{equation}
so that $\nu_p(t)=\nu_p(i)=n$ and $p\nmid t'$, where $\nu_p$ denotes the $p$-adic valuation.

The key objects are the modular curves $X_0(\ell)$, whose non-cuspidal points over a field classify pairs $(E,C)$ of an elliptic curve with a cyclic subgroup of order $\ell$. Such a point is called \emph{supersingular} if the underlying curve is. Throughout, $\Fpb$ denotes an algebraic closure of $\F_p$ and $\Fpp\subset\Fpb$ the field with $p^2$ elements. For a prime $\ell\ne p$, let $V_\ell$ denote the classical operator on $q$-expansions $f(q)\mapsto f(q^\ell)$, which maps forms of level one to forms on $\Gamma_0(\ell)$ of the same weight, and set
\begin{equation}
R_\ell\ :=\ \frac{V_\ell(\Delta)}{\Delta},
\end{equation}
a modular \emph{function} on $X_0(\ell)$. In Section~\ref{sec:containment}, we recall its moduli interpretation: at a point of $X_0(\ell)$ corresponding to a pair $(E,C)$, the value of $R_\ell$ is the normalized ratio of discriminants along the $\ell$-isogeny $E\to E/C$.

\begin{theorem}[Supersingular containment]\label{thm:containment}
Let $p\ge 5$ and let $\ell\ne p$ be any prime. Then for every supersingular point $\wt P$ of $X_0(\ell)_{\Fpb}$, we have
\[
R_\ell(\wt P)^{\,m}\ =\ 1,
\]
where $m=(p^2-1)/12.$
\end{theorem}

Theorem~\ref{thm:containment} follows from a rationality theorem for supersingular discriminants, which is of independent interest. For an elliptic curve $E$ over a field with a nonzero invariant differential $\omega$, we write $\Delta(E,\omega)$ for the value of the usual weight-$12$ discriminant modular form at $(E,\omega)$ (see \S\ref{subsec:hasse} for the precise conventions). It is nonzero, and it satisfies $\Delta(E,\lam\omega)=\lam^{-12}\Delta(E,\omega)$, so the quantity $\kappa(E)$ below does not depend on the choice of $\omega$. We write $[n]$ for the multiplication-by-$n$ endomorphism. Recall that a supersingular elliptic curve $E$ over $\Fpp$ has $p^2$-power Frobenius endomorphism $\pi_E\in\End(E)$ satisfying 
$$\pi_E^2-a\pi_E+p^2=0,$$
 with $|a|\le 2p$. The extreme case $a=-2p$  (i.e.\ $\pi_E=[-p]$) occurs for a suitable model of \emph{every} supersingular $j$-invariant (see Lemma~\ref{lem:deuring}). We call such models \emph{Deuring models}.

\begin{theorem}[Rationality of supersingular discriminants]\label{thm:kappa}
Let $p\ge 5$, let $E/\Fpp$ be a supersingular elliptic curve with $\pi_E=[-p]$, and let $\omega$ be any nonzero $\Fpp$-rational invariant differential on $E$. Then we have that
\[
\kappa(E):=\Delta(E,\omega)^{(p^2-1)/12}\ =\ (-1)^{(p+1)/2}\ =\
\begin{cases}
+1 & \text{if } p\equiv 3\pmod 4,\\
-1 & \text{if } p\equiv 1\pmod 4.
\end{cases}
\]
In particular $\kappa(E)$ depends on neither $E$ nor $\omega$.
\end{theorem}

\subsection*{Overview of the proofs}
Let $L_\ell:=V_\ell-\iota$, with $\iota$ the inclusion from level one to level $\Gamma_0(\ell)$. For $f=\Delta^t$, condition (C4) asks that $A_p^{\,j}\mid L_\ell(\Delta^t)$ modulo $p$ for every $\ell\ne p$.

At level $2$, Larson's generators satisfy $\mu:=\delta^2-\eps$, $\Delta=64\mu\eps^2$, and $V_2(\Delta)=\mu^2\eps$, so we have
\[
L_2(\Delta^t)
=\mu^t\eps^t\bigl(\mu^t-64^t\eps^t\bigr)
=\mu^t\eps^t\!\!\prod_{\rho^{\,t'}=64^{t'}}\!\!(\mu-\rho\eps)^{p^n}
\]
in $\Fpb[\delta,\eps]$. Igusa's factorization of $A_p$ on $X_0(2)$ reduces the positive divisibility to the supersingular identity $R_2(\wt P)^m=1$; the negative direction follows from the exact multiplicities of the same factors.

Theorem~\ref{thm:containment} proves the identity $R_\ell(\wt P)^m=1$ for all $\ell\ne p$. After replacing $(E,C)$ by a Deuring model, the isogeny $E\to E/C$ is defined over $\Fpp$, and the discriminant weights give
\[
R_\ell(\wt P)^m=\frac{\kappa(E/C)}{\kappa(E)}.
\]
Theorem~\ref{thm:kappa}, proved from the level $3$ and level $4$ eta powers whose powers recover $\Delta$, makes this ratio equal to $1$.

Thus Section~\ref{sec:behrens} recalls the congruence conditions, Section~\ref{sec:levels} supplies the level $2$ calculation, Section~\ref{sec:containment} proves the supersingular discriminant and containment theorems, Section~\ref{sec:proof} deduces Theorem~\ref{thm:main}, and Section~\ref{sec:examples} offers examples when $p=11$ and $13$.

\section*{Acknowledgements}
\noindent The author thanks Donald Larson for bringing this problem to his attention, and he thanks Ashvin Swaminathan and Nick Kuhn for comments on early drafts of this paper. The author also thanks the Thomas Jefferson Fund, the NSF (DMS-2002265 and DMS-2055118) and the Simons Foundation (SFI-MPS-TSM-00013279) for their generous support.

\section{Behrens' conditions and modular forms modulo $p$}\label{sec:behrens}

This section fixes the congruence language used throughout the proof. We recall only the parts of Behrens' construction that are needed to test the candidate form $\Delta^t$.

\subsection{Behrens' conditions}\label{subsec:conditions}
We begin by recording the indexing convention and the four conditions that characterize the Behrens forms in Larson's formulation. This will allow the rest of the paper to be purely arithmetic.

Fix a prime $p\ge 5$. For an integer $N\ge1$ let $\Gamma(N)\subset\SL_2(\Z)$ denote the subgroup of matrices congruent to the identity modulo $N$, and for a prime $\ell$ let $\Gamma_0(\ell)$ denote the subgroup of matrices that are upper triangular modulo $\ell$. For a power series $f=\sum_n a_nq^n$ we write $\ord_q(f)$ for the smallest $n$ with $a_n\ne0$. For a congruence subgroup $\Gamma$ we write $M_k(\Gamma,R)$ for the space of holomorphic modular forms of weight $k$ with coefficients in the ring $R$, and $M_*(\Gamma,R)=\bigoplus_k M_k(\Gamma,R)$; bars denote reduction modulo $p$. For a prime $\ell$, the operator $V_\ell\colon f(q)\mapsto f(q^\ell)$ maps $M_k(\SL_2(\Z),\Z)$ into $M_k(\Gamma_0(\ell),\Z)$ (see, e.g., \cite{OnoCBMS}), and we set
\begin{equation}\label{eq:L}
L_\ell\ :=\ V_\ell-\iota,
\end{equation}
where $\iota\colon M_k(\SL_2(\Z),\Z)\hookrightarrow M_k(\Gamma_0(\ell),\Z)$ is the inclusion. This is Larson's operator \cite{Larson}. Note that $V_\ell$ is a ring homomorphism on $q$-expansions.

For the pairs of integers $(i,j)$ enumerated by Larson \cite[\S2, Lem.~2.1]{Larson} --- the indices of the order $p$ elements $\beta_{i/j}$ of the divided beta family --- Behrens' theorem \cite[Thm.~1.3]{BehrensCongruences} (in the formulation of \cite[\S2]{Larson}) attaches to $\beta_{i/j}$ a modular form $f_{i/j}$ modulo $p$ of weight $12t$, $t=i(p^2-1)/12$ (the reduction $\bar f_{i/j}$ of an integral form). We recall the enumeration using the low-index convention forced by Larson's examples. Put
\[
b_m:=
\begin{cases}
0,&m<0,\\
1,&m=0,\\
p^m+p^{m-1}-1,&m\ge1,
\end{cases}
\]
and write $i=rp^n$ with $(r,p)=1$. Then $\beta_{i/j}$ is a divided beta element of order $p$ --- and we call the pair $(i,j)$ \emph{admissible} --- exactly for the following $j$ \textup{\cite[Lem.~2.1]{Larson}; cf.\ \cite[Ex.~2.2]{Larson} for $n=0$ and \cite[Ex.~2.3]{Larson} for the exclusion of $j=p$ when $n=2$}:
\begin{itemize}[leftmargin=2.1em]
\item if $r=1$: for $1\le j\le p^{n}$ with $j\notin\{p,2p,\dots,b_{n-2}p\}$;
\item if $r>1$: for $1\le j\le b_{n}$ with $j\notin\{p,2p,\dots,b_{n-2}p\}$.
\end{itemize}
Therefore,  the excluded set is empty for $n\le1$. For $n=2,$ it consists of the single index $p$. In all cases the excluded indices are $\le b_{n-2}p<p^{n}$, so they lie strictly inside the range $j\le p^n$; the region $j>p^n$ contains no exclusions. The convention $b_0=1$ records Larson's standard elements $\beta_i=\beta_{i/1}$ for all $i\ge1$. Each admissible $f_{i/j}$ is characterized by the following conditions, and is determined by them only up to multiplication by a unit modulo $p$ and addition of a form of the same weight and $q$-order that reduces to $0$ modulo $p$ \cite[\S1]{Larson}:
\begin{description}[leftmargin=2.1em]
\item[(C1)] $\bar f_{i/j}\ne 0$.
\item[(C2)] The integer $12\ord_q(\bar f_{i/j})$ is either strictly greater than $12t-(p-1)j$, or is equal to $12t-(p-1)j-2$; see \cite[Thm.~1.3]{BehrensCongruences} and \cite[Lem.~4.3]{Larson}.
\item[(C3)] No form of weight smaller than $12t$ is congruent to $f_{i/j}$ modulo $p$.
\item[(C4)] For every prime $\ell\ne p$,
\[
A_p^{\,j}\ \Bigm|\ L_\ell(\bar f_{i/j})\qquad\text{in } M_\ast(\ell)_{\Z/p},
\]
where $A_p$ is the Hasse invariant of \S\ref{subsec:hasse} and $M_\ast(\ell)_{\Z/p}$ denotes the graded ring of modular forms modulo $p$ on $\Gamma_0(\ell)$ in the geometric sense of Katz \cite{Katz350}. (Remark~\ref{rem:rings} compares this ring with the ring of reductions of integral forms; at $\ell=2$, the level at which clause (ii) of Theorem~\ref{thm:main} is proved, both readings agree with Larson's ring $\F_p[\delta,\eps]$.)
\end{description}
Any form satisfying (C1)--(C4) may serve as $f_{i/j}$; conversely, a form failing any of the conditions cannot equal $f_{i/j}$. We refer to \cite{BehrensCongruences,Larson} for the precise indexing and for the topological content of these conditions.

\subsection{The Hasse invariant}\label{subsec:hasse}
The Hasse invariant is the geometric divisor that detects the supersingular locus. We also fix the evaluation convention for the discriminant, since later arguments compare discriminants along isogenies.

For $p\ge5$ the Eisenstein series $E_{p-1}$ has $p$-integral $q$-expansion with $E_{p-1}\equiv 1\pmod p$, by the von Staudt--Clausen theorem (see \cite{SwD}). We write
\[
A_p\ :=\ \overline{E}_{p-1}\in M_{p-1}(\SL_2(\Z),\F_p)
\]
for its reduction, the \emph{Hasse invariant}. Following Katz \cite[\S1]{Katz350}, a modular form of weight $k$ and level $\Gamma$ over a ring $R_0$ is a rule assigning to each triple $(E/R,\ \text{level structure},\ \omega)$ --- with $R$ an $R_0$-algebra, $E/R$ an elliptic curve, and $\omega$ a nowhere-vanishing invariant differential --- an element of $R$, compatibly with base change and isomorphisms and satisfying $f(E,\cdot,\lam\omega)=\lam^{-k}f(E,\cdot,\omega)$ for $\lam\in R^\times$. Classical forms with $R_0$-integral $q$-expansions define such rules, and evaluation on the Tate curve with its canonical differential recovers the $q$-expansion \cite[\S1.1--1.7]{Katz350}, \cite{KatzMazur}. The basic example is the discriminant. For a short Weierstrass equation
\[
E\colon y^2=x^3+Ax+B
\]
with its standard invariant differential $\omega=dx/(2y)$, one has
\[
\Delta(E,\omega):=-16\,(4A^3+27B^2).
\]
Equivalently, if the same equation is paired with the differential $dx/y=2\omega$, the modular-form value is $2^{-12}\cdot[-16(4A^3+27B^2)]$. General Weierstrass models are covered by the standard formulas \cite[\S III.1]{SilvermanAEC}, and rescaling the differential gives $\Delta(E,\lam\omega)=\lam^{-12}\Delta(E,\omega)$. The $q$-expansion of this rule --- its value on the Tate curve with its canonical differential --- is the series $\Delta(q)$ \cite[Ch.~8]{KatzMazur}. We will use the following standard facts.

\begin{lemma}\label{lem:hasse}
Let $p\ge 5$ and let $k$ be a field of characteristic $p$. Then the following are true.
\begin{enumerate}
\item[(a)] For an elliptic curve $E/k$ with invariant differential $\omega$, one has $A_p(E,\omega)=0$ if and only if $E$ is supersingular \cite[\S2]{Katz350}, \cite[Ch.~12]{KatzMazur}.
\item[(b)] For every prime $\ell$ and every cusp of $X_0(\ell)$, the constant term of the expansion of $A_p$ at that cusp equals $1$. In particular $A_p$ does not vanish at any cusp.
\item[(c)] If $E\colon y^2=f(x)$ with $f$ a separable cubic, then $A_p(E,dx/y)$ equals the coefficient of $x^{p-1}$ in $f(x)^{(p-1)/2}$ \cite[\S12.4]{KatzMazur}, \cite[\S V.4]{SilvermanAEC}. The same value is obtained for $dx/(2y)$, since $A_p$ has weight $p-1$ and $2^{p-1}=1$ in characteristic $p$.
\item[(d)] For the Legendre curve $E_\lam\colon y^2=x(x-1)(x-\lam)$,
\[
A_p\bigl(E_\lam,\tfrac{dx}{2y}\bigr)\;=\;A_p\bigl(E_\lam,\tfrac{dx}{y}\bigr)\;=\;(-1)^{(p-1)/2}H_p(\lam),\qquad
H_p(\lam):=\sum_{i=0}^{(p-1)/2}\binom{(p-1)/2}{i}^{2}\lam^{i}\in\F_p[\lam].
\]
\item[(e)] $H_p$ is separable, and its roots in $\Fpb$ are exactly the supersingular values of $\lam$ \cite[Thm.~V.4.1]{SilvermanAEC}; the separability is a theorem of Igusa \cite{Igusa} (cf.\ \cite{BrillhartMorton,Deuring}).
\end{enumerate}
\end{lemma}

\begin{proof}
Parts (a), (c) and (e) carry the indicated references; in (a) and (c) we also use that the geometric Hasse invariant --- defined via the action of absolute Frobenius on $H^1(E,\mathcal O_E)$, as in \cite[\S12.4]{KatzMazur} --- coincides with $\overline{E}_{p-1}$: both are weight $p-1$ forms of level one over $\F_p$ with $q$-expansion $1$ (for the geometric Hasse invariant this is the Tate-curve computation in \cite[\S12.4]{KatzMazur}), and the $q$-expansion map is injective in each weight. For (b): if $\gamma\in\SL_2(\Z)$ carries $\infty$ to the given cusp, the expansion of $A_p$ there is computed from $E_{p-1}|_{p-1}\gamma=E_{p-1}$, whose constant term is $1$. For (d), write $w:=\tfrac{p-1}{2}$. By (c) we must extract the coefficient of $x^{p-1}=x^{2w}$ in $x^{w}(x-1)^{w}(x-\lam)^{w}$, i.e.\ the coefficient of $x^{w}$ in $(x-1)^w(x-\lam)^w$, which equals
\[
\sum_{a+b=w}\binom{w}{a}(-1)^{w-a}\binom{w}{b}(-\lam)^{w-b}
=(-1)^{w}\sum_{a=0}^{w}\binom{w}{a}^{2}\lam^{a}=(-1)^{w}H_p(\lam). \qedhere
\]
\end{proof}

\subsection{Two lemmas for (C2) and (C3)}
The following elementary facts dispose of the first three Behrens conditions for the candidate $\Delta^t$. They are independent of the level-$\ell$ divisibility condition.

\begin{lemma}\label{lem:maxord}
Let $g\in M_{12t}(\SL_2(\Z),\Z)$ with $\bar g\ne 0$. Then $\ord_q(\bar g)\le t$. In particular $\Delta^t$, whose reduction has $\ord_q=t$, attains the maximal possible $q$-order in its weight.
\end{lemma}

\begin{proof}
Suppose $\ord_q(\bar g)\ge 1$, i.e.\ the constant term $c$ of $g$ is divisible by $p$. Then $g_1:=g-c\,E_4^{3t}$ lies in $M_{12t}(\SL_2(\Z),\Z)$ and has vanishing constant term, hence is a cusp form, hence $g_1=\Delta h$ with $h\in M_{12t-12}(\SL_2(\Z),\Z)$: the quotient has integral $q$-expansion because $\Delta=q\cdot(\text{unit in }\Z[[q]])$, and it is a holomorphic form of weight $12t-12$. Reducing, $\bar g=\bar\Delta\bar h$ and $\ord_q(\bar h)=\ord_q(\bar g)-1$. If $\ord_q(\bar g)\ge t+1$, iterating $t+1$ times produces a form of weight $12t-12(t+1)=-12<0$ with nonzero reduction, which is impossible.
\end{proof}

\begin{lemma}\label{lem:filtration}
There is no $g\in M_{k'}(\SL_2(\Z),\Z)$ with $k'<12t$ and $\bar g=\overline{\Delta^{t}}$. Hence $\Delta^t$ satisfies \textup{(C3)}.
\end{lemma}

\begin{proof}
Suppose such a $g$ exists. By the theory of Serre and Swinnerton-Dyer \cite{SwD} (see also \cite[\S4]{Katz350}), two level-one forms with equal nonzero reductions have weights congruent modulo $p-1$; write $e=(12t-k')/(p-1)\ge1$. Then $A_p^{\,e}\,\bar g$ has weight $12t$ and the same $q$-expansion as $\overline{\Delta^t}$, so $\overline{\Delta^t}=A_p^{\,e}\,\bar g$ by injectivity of $q$-expansions in each weight; in particular $A_p\mid\overline{\Delta^t}$ in the graded ring $M_*(\SL_2(\Z),\Z)\otimes\F_p$. Evaluating at a supersingular pair $(E,\omega)$ --- these exist since $\deg H_p=(p-1)/2\ge1$, by Lemma~\ref{lem:hasse}(e) --- the right side vanishes by Lemma~\ref{lem:hasse}(a), while $\Delta(E,\omega)^t\ne0$. This contradiction proves the lemma.
\end{proof}

\section{Level 2: Connection to the Hasse invariant}\label{sec:levels}

This section carries out the explicit calculation at level $2$. The point is to express both $L_2(\Delta^t)$ and the pullback of $A_p$ in Larson's two generators.

\subsection{Larson's generators and the coordinate $x$}\label{subsec:larsonring}
We first recall the weighted polynomial presentation of the level-$2$ ring. The quotient $x=\mu/\eps$ will serve as the coordinate in which supersingular factors are read off.

Following Larson \cite[\S3]{Larson}, set
\[
\delta:=\frac{2E_2(2\tau)-E_2(\tau)}{4}\in M_2(\Gamma_0(2)),\qquad
\eps:=\frac{1}{16}\,\frac{\eta(\tau)^{16}}{\eta(2\tau)^{8}}\in M_4(\Gamma_0(2)),\qquad \mu:=\delta^2-\eps,
\]
where $E_2(\tau)=1-24\sum_{n\ge1}\sigma_1(n)q^n$ is the quasimodular Eisenstein series of weight $2$ and $\eta(\tau)=q^{1/24}\prod_{n\ge1}(1-q^n)$ is the Dedekind eta function. The $q$-expansions are
\[
\delta=\tfrac14+6q+6q^2+24q^3+\cdots,\qquad
\eps=\tfrac1{16}-q+7q^2-28q^3+\cdots,\qquad
\mu=4q+32q^2+112q^3+\cdots .
\]
Then $M_*(\Gamma_0(2),\Z[1/2])=\Z[1/2][\delta,\eps]$ is a weighted polynomial ring \cite[\S3]{Larson}, so for odd $p$,
\begin{equation}\label{eq:ring}
M_*(2)_{\Z/p}\ :=\ M_*(\Gamma_0(2),\Z[1/2])\otimes\F_p\ =\ \F_p[\delta,\eps]
\end{equation}
(a weighted polynomial ring in which $\delta$ has weight $2$ and $\eps$ has weight $4$), and reduction modulo $p$ is computed on $q$-expansions.

The two identities below are the bridge between Larson's level-$2$ ring and the operator $L_2$.

\begin{lemma}\label{lem:identities}
In $\Z[1/2][\delta,\eps]$ one has the exact identities
\begin{equation}\label{eq:DeltaVDelta}
\Delta\ =\ 64\,\mu\,\eps^{2},\qquad V_2(\Delta)\ =\ \mu^{2}\eps .
\end{equation}
\end{lemma}

\begin{proof}
All four forms lie in $M_{12}(\Gamma_0(2),\Z[1/2])$. By the Sturm bound \cite{Sturm}, an element of $M_{12}(\Gamma_0(2))$ vanishing to order greater than $\tfrac{12}{12}\,[\SL_2(\Z):\Gamma_0(2)]=3$ at $\infty$ is zero, so it suffices to compare $q$-expansions through $q^3$ --- a finite computation from the definitions above (cf.\ \cite[\S3--4]{Larson}); one finds $64\mu\eps^2=q-24q^2+252q^3+\cdots=\Delta$ and $\mu^2\eps=q^2+O(q^4)=\Delta(q^2)$ through the required order.
\end{proof}

The next lemma records the geometry of the coordinate $x$ and identifies the exceptional elliptic point.

\begin{lemma}\label{lem:hauptmodul}
Work over $\Fpb$ ($p\ge5$), and set
\[
R_2:=\frac{V_2(\Delta)}{\Delta}=\frac{\mu}{64\,\eps},\qquad x:=64R_2=\frac{\mu}{\eps}.
\]
Then the following are true.
\begin{enumerate}
\item[(a)] We have that
$$R_2=q+24q^2+300q^3+\cdots.
$$
 It is holomorphic and nonvanishing away from the cusps of $X_0(2)$, and its divisor is $(\infty)-(0)$, where $\infty$ and $0$ denote the two cusps.
\item[(b)] We have that $x$ defines an isomorphism $X_0(2)_{\Fpb}\xrightarrow{\ \sim\ }\PP^1$, with $x(\infty)=0$ and $x(0)=\infty$. In particular, for every $\xi\in\Fpb$ the function $x-\xi$ is a uniformizer at the unique point where $x=\xi$.
\item[(c)] We have that $\Delta$, $\mu$ and $\eps$ are nonvanishing at every non-cuspidal point. Moreover, $\mu$ vanishes at the cusp $\infty$, while $\eps$ vanishes at the cusp $0$ and satisfies $\eps(\infty)=\tfrac1{16}\ne0$. Furthermore, we have that
\[
x+1=\frac{\delta^2}{\eps},
\]
and the unique point with $x=-1$ is the unique point of $X_0(2)$ at which $\delta$ vanishes.
\end{enumerate}
\end{lemma}

\begin{proof}
(a) The $q$-expansion is that of $\Delta(q^2)/\Delta(q)$. On the open modular curve both $\Delta$ and $V_2\Delta$ are values of $\Delta$ on elliptic curves (Lemma~\ref{lem:Vell} below), hence nonvanishing; so the divisor of the function $R_2$ is supported on the two cusps. The cusp $\infty$ of $X_0(2)$ has width one, so $q$ is a formal parameter there \cite[Ch.~8--10]{KatzMazur}, and $\ord_\infty(R_2)=1$. A principal divisor has degree zero, so $\operatorname{div}(R_2)=(\infty)-(0)$.

(b) By (a), $x=64R_2$ is a nonconstant map $X_0(2)_{\Fpb}\to\PP^1$ with a single simple pole; hence it has degree one and is an isomorphism, with zero at the cusp $\infty$ and pole at the cusp $0$.

(c) $\Delta$ is nonvanishing off the cusps, and $\Delta=64\mu\eps^2$ forces the same for $\mu$ and $\eps$. The displayed $q$-expansions give $\mu(\infty)=0$ and $\eps(\infty)=\tfrac1{16}\ne0$ directly. At the cusp $0$, the function $x=\mu/\eps$ has a simple pole by (b), so $\ord_0(\eps)=\ord_0(\mu)+1\ge1$, i.e.\ $\eps$ vanishes there. The identity $x+1=(\mu+\eps)/\eps=\delta^2/\eps$ is immediate from $\mu=\delta^2-\eps$. If $\delta(y)=0$ at a point $y$, then $y$ is not a cusp: the constant term of $\delta$ at $\infty$ is $\tfrac14$, and at the cusp $0$ it is again nonzero, since $2E_2(2\tau)-E_2(\tau)$ is sent to $-2z^2\bigl(2E_2(2z)-E_2(z)\bigr)$ under $\tau\mapsto-1/(2z)$ (the two nonmodular error terms of $E_2$ cancel), so the expansion of $\delta$ at $0$ is a nonzero multiple of its expansion at $\infty$. Hence $\eps(y)\ne0$ and $x(y)+1=\delta(y)^2/\eps(y)=0$; conversely at the point $x=-1$ one gets $\delta^2=\eps\cdot(x+1)=0$.
\end{proof}

\subsection{The double cover by $X(2)$ and the factorization of $A_p$}\label{subsec:factorAp}
We now pass to full level $2$, where the Hasse invariant is the classical Deuring polynomial in the Legendre parameter. Pushing this information back to $X_0(2)$ gives the needed factorization in $x$.

Let $X(2)$ denote the modular curve of full level $2$ structures over $\Z[1/2]$: it is $\PP^1$ with affine coordinate $\lam$, the universal curve away from the cusps $\lam\in\{0,1,\infty\}$ being the Legendre curve $E_\lam$ with the ordered basis $P_1=(0,0)$, $P_2=(1,0)$ of $E_\lam[2]$; the coordinate $\lam$ classifies such triples $(E;P_1,P_2)$ up to isomorphism \cite{DeligneRapoport,KatzMazur}. Let
\[
\pi\colon X(2)\longrightarrow X_0(2),\qquad (E;P_1,P_2)\longmapsto (E,\langle P_1\rangle).
\]

The covering map is unramified away from the familiar order-$4$ elliptic point, and this is exactly what accounts for the possible squared factor in $A_p$.

\begin{lemma}\label{lem:cover}
Over $\Fpb$ ($p\ge5$), $\pi$ is a finite separable morphism of degree $2$ of smooth projective curves. For a non-cuspidal point $y=(E,C)$ of $X_0(2)(\Fpb)$, the fiber $\pi^{-1}(y)$ consists of two distinct points, at each of which $\pi$ is unramified, unless $j(E)=1728$ and $C$ is the subgroup of $E[2]$ fixed by an automorphism of order $4$; there is exactly one such point $y_2$ (the elliptic point of $X_0(2)$), its fiber is a single point at which the ramification index is $2$, and
\[
x(y_2)=-1 .
\]
\end{lemma}

\begin{proof}
Non-cuspidal $\Fpb$-points of the coarse curves classify isomorphism classes \cite{KatzMazur}. Given $(E,C)$, the triples above it are $(E;P_1,P_2)$ with $P_1$ the nonzero point of $C$ and $P_2$ one of the two points of $E[2]\smallsetminus C$; two such triples define the same point of $X(2)$ if and only if some $\alpha\in\Aut(E)$ fixes $P_1$ and exchanges the two candidates for $P_2$. If $j\ne 0,1728$ then $\Aut(E)=\{\pm1\}$ acts trivially on $E[2]$, so the fiber has two points. If $j=0$, an automorphism $\rho$ of order $3$ cannot act trivially on $E[2]$ (else $0=\rho^2+\rho+1$ would act as multiplication by $3=1\ne0$ on $E[2]$), so $\rho$ permutes the three nonzero $2$-torsion points cyclically and fixes none of them; hence $\Aut(E,C)=\{\pm1\}$ and again the fiber has two points. If $j=1728$, write $E\colon y^2=x^3+ax$ and $[i](x,y)=(-x,iy)$: then $[i]$ fixes $(0,0)$ and exchanges the other two nonzero $2$-torsion points. Thus for $C=\langle(0,0)\rangle$ the two triples are isomorphic and the fiber is one point, while the other two choices of $C$ are exchanged by $[i]$, hence give a single point of $X_0(2)$ whose fiber has two points. The point $y_2=(E,\langle(0,0)\rangle)$ is unique up to isomorphism. Since $\deg\pi=2$ is prime to $p$, $\pi$ is separable, so $\sum_{Q\mapsto y}e_Q=2$ for every $y$, where $e_Q$ denotes the ramification index of $\pi$ at $Q$; the fiber counts now give the ramification statement.

It remains to prove $x(y_2)=-1$. The forms $\delta,\eps$ are algebraic modular forms over $\Z[1/2]$ (classical forms of weight $\ge2$ with $\Z[1/2]$-integral $q$-expansions; \cite[\S1.7]{Katz350}), so
\[
h(\lam):=\delta\bigl(E_\lam,\langle P_1\rangle,\tfrac{dx}{2y}\bigr),\qquad
e(\lam):=\eps\bigl(E_\lam,\langle P_1\rangle,\tfrac{dx}{2y}\bigr)
\]
are regular functions on the $\lam$-line minus the cusps, and $e(\lam)$ is nonvanishing there by Lemma~\ref{lem:hauptmodul}(c), while
\[
(x\circ\pi)(\lam)+1=\frac{h(\lam)^2}{e(\lam)} .
\]
Hence, the function $x\circ\pi+1$ has \emph{even} order of vanishing at every non-cuspidal point of $X(2)$. The value $x=-1$ is not a cusp (the cusps are $x\in\{0,\infty\}$ by Lemma~\ref{lem:hauptmodul}), so there is a non-cuspidal $Q\in X(2)(\Fpb)$ with $\pi(Q)$ equal to the point $x=-1$; then $\ord_Q(x\circ\pi+1)=e_Q\cdot\ord_{x=-1}(x+1)=e_Q$ is even, so $e_Q=2$ and $Q$ is ramified. By the first part, $\pi(Q)=y_2$, i.e.\ $x(y_2)=-1$.
\end{proof}

\begin{proposition}\label{prop:Apfactor}
Let $p\ge5$ and work in the weighted polynomial ring $\Fpb[\delta,\eps]$. Put $\alpha_p:=0$ if $p\equiv1\pmod4$ and $\alpha_p:=1$ if $p\equiv3\pmod4$, and $s:=\dfrac{p-1-2\alpha_p}{4}\in\Z_{\ge1}$. Then there exist a constant $c\in\F_p^\times$ and \emph{distinct} elements $\xi_1,\dots,\xi_s\in\Fpb\smallsetminus\{0,-1\}$ such that
\[
A_p\ =\ c\;\delta^{\,\alpha_p}\prod_{i=1}^{s}\bigl(\mu-\xi_i\,\eps\bigr).
\]
Moreover the supersingular points of $X_0(2)_{\Fpb}$ are precisely the $s$ points with $x=\xi_i$, together with the point $x=-1$ when $\alpha_p=1$.
\end{proposition}

\begin{proof}
Every nonzero homogeneous element of $\Fpb[\delta,\eps]$ (weights $2$ and $4$) factors as a constant times $\delta^{a}$ times a product of weight-$4$ linear forms in $(\mu,\eps)$: writing it as $\delta^{a}g(\delta^2,\eps)$ with $g$ a binary form in $(\delta^2,\eps)$ and factoring $g$ into linear forms over $\Fpb$, each factor $\alpha\delta^2+\beta\eps=\alpha\mu+(\alpha+\beta)\eps$ is linear in $(\mu,\eps)$, and factors proportional to $\mu+\eps=\delta^2$ are absorbed into $\delta^{a}$. Applying this to $A_p$, which lies in $\F_p[\delta,\eps]$ by \eqref{eq:ring}, gives
\[
A_p=c\,\delta^{a}\,\mu^{m_0}\,\eps^{m_\infty}\prod_{i}(\mu-\xi_i\eps)^{m_i},\qquad \xi_i\in\Fpb\smallsetminus\{0,-1\}\ \text{distinct},\ c\in\F_p^\times .
\]
The constant term of the $q$-expansion of $A_p$ at $\infty$ is $1$ while that of $\mu$ is $0$; hence $m_0=0$. If $m_\infty\ge1$ then $A_p$ would vanish at the cusp $0$, where $\eps$ vanishes (Lemma~\ref{lem:hauptmodul}(c)), contradicting Lemma~\ref{lem:hasse}(b); hence $m_\infty=0$. By Lemma~\ref{lem:hasse}(a) the non-cuspidal zero locus of $A_p$ on $X_0(2)$ is exactly the supersingular locus; by Lemma~\ref{lem:hauptmodul} the non-cuspidal zero locus of the factor $\mu-\xi_i\eps$ is the single point $x=\xi_i$ (at non-cuspidal points $\eps\ne0$, so $\mu-\xi_i\eps=0$ if and only if $x=\xi_i$), and the zero locus of $\delta$ is the point $x=-1$. Hence the supersingular points are the points $x=\xi_i$ together with $x=-1$ exactly when $a\ge1$, and it remains to show $m_i=1$ for all $i$ and $a\le 1$.

Consider the weight-zero ratio $F:=A_p^{12}/\Delta^{p-1}$, a rational function on $X_0(2)_{\Fpb}$. Using $\mu-\xi\eps=\eps(x-\xi)$, $\delta^2=\eps(x+1)$ and $\Delta=64\mu\eps^2=64\eps^3 x$, and noting that the total exponent of $\eps$ is $6a+12\sum_i m_i-3(p-1)=3\bigl(2a+4\sum_i m_i-(p-1)\bigr)=0$ by comparing weights, we get
\[
F\ =\ \frac{c^{12}}{64^{\,p-1}}\;(x+1)^{6a}\;\prod_i (x-\xi_i)^{12m_i}\;x^{-(p-1)} .
\]
Since $x$ is a global coordinate (Lemma~\ref{lem:hauptmodul}(b)),
\[
\ord_{x=\xi_i}(F)=12m_i,\qquad \ord_{x=-1}(F)=6a .
\]
On the other hand, $A_p$ and $\Delta$ are forms of level one, so the value of $F$ at a point does not depend on the level structure, and the pullback of $F$ along $\pi$ is computed on the Legendre family: by Lemma~\ref{lem:hasse}(d) and $\Delta(E_\lam,dx/(2y))=16\lam^2(\lam-1)^2$,
\[
(\pi^*F)(\lam)\ =\ \frac{H_p(\lam)^{12}}{\bigl(16\lam^2(\lam-1)^2\bigr)^{p-1}} .
\]
At a supersingular $\lam_0$ (which is not a cusp), $H_p$ has a simple zero by Lemma~\ref{lem:hasse}(e), so $\ord_{\lam_0}(\pi^*F)=12$. Since $\ord_{\lam_0}(\pi^*F)=e_{\lam_0}\cdot\ord_{\pi(\lam_0)}(F)$, Lemma~\ref{lem:cover} gives: at a supersingular point $x=\xi_i$ (non-elliptic, $e=1$) we get $12m_i=12$, i.e.\ $m_i=1$; and if $a\ge1$, so that the point $x=-1$ is supersingular, its unique preimage is the ramification point ($e=2$) and $12=2\cdot 6a$, i.e.\ $a=1$. Finally the weight count $2a+4s=p-1$ forces $a\equiv\frac{p-1}{2}\pmod 2$, whence $a=\alpha_p$ as defined and $s=(p-1-2\alpha_p)/4$, which is a positive integer for every $p\ge5$.
\end{proof}

\begin{remark}
Proposition~\ref{prop:Apfactor} recovers, up to units, Larson's factorizations $A_{11}=\delta(\mu+3\eps)(\mu+4\eps)$ and $A_{13}=c\,(\mu+12\eps)(\mu^2+5\mu\eps+\eps^2)$ \cite[\S5]{Larson}; see Section~\ref{sec:examples}. It also shows that the elliptic point of $X_0(2)$ is supersingular exactly when $p\equiv 3\pmod 4$, in accordance with the classical criterion for $j=1728$.
\end{remark}

\section{Supersingular values of the isogeny ratio}\label{sec:containment}

This section proves the supersingular identity that drives the divisibility in all levels. The argument is geometric: the isogeny ratio is compared with a field-independent power of the discriminant.

\subsection{The moduli interpretation}
We first interpret $R_\ell$ as a ratio of discriminants along the universal $\ell$-isogeny. This is the form in which supersingular rationality can be applied.

Fix primes $p\ge5$ and $\ell\ne p$. For a field $k$ of characteristic $p$, a pair $(E,C)$ consisting of an elliptic curve over $k$ and a subgroup $C$ of order $\ell$, and a nonzero invariant differential $\omega$, let $\varphi\colon E\to E':=E/C$ be the quotient isogeny --- which exists over $k$, with explicit equations, by V\'elu \cite{Velu}; cf.\ \cite[\S III.4]{SilvermanAEC} --- and let $\omega'$ be the unique invariant differential on $E'$ with $\varphi^*\omega'=\omega$. (It exists and is nonzero because $\deg\varphi=\ell$ is prime to $p$, so $\varphi$ is separable.) Define
\[
\Delta^{[\ell]}(E,C,\omega)\ :=\ \ell^{-12}\,\Delta(E',\omega') .
\]

\begin{lemma}\label{lem:Vell}
$\Delta^{[\ell]}$ is a modular form of weight $12$ and level $\Gamma_0(\ell)$ over $\F_p$ in the sense of Katz, with $q$-expansion $\Delta(q^\ell)$; it is the reduction of the classical $V_\ell(\Delta)$, and
\[
R_\ell:=\frac{\Delta^{[\ell]}}{\Delta}
\]
is the rational function on $X_0(\ell)_{\Fpb}$ whose expansion at the cusp $\infty$ is $V_\ell(\Delta)/\Delta$. For $\ell=2$ one has $R_2=\mu/(64\eps)=x/64$ in the notation of Section~\textup{\ref{sec:levels}}.
\end{lemma}

\begin{proof}
The rule $\Delta^{[\ell]}$ is compatible with base change and isomorphisms because the formation of $E/C$ and of $\omega'$ is; and replacing $\omega$ by $\lam\omega$ replaces $\omega'$ by $\lam\omega'$ and multiplies the value by $\lam^{-12}$, giving weight $12$. To compute the $q$-expansion, evaluate on the Tate curve: over $\F_p((q))$ one has $\mathrm{Tate}(q)=\mathbb G_m/q^{\Z}$ with canonical differential $\omega_{\mathrm{can}}=du/u$ and canonical subgroup $\mu_\ell$ of $\ell$-th roots of unity (the group scheme $\mu_\ell$, not to be confused with the level-$2$ form $\mu$ of Section~\ref{sec:levels}), the cusp $\infty$ of $X_0(\ell)$ corresponds to $(\mathrm{Tate}(q),\mu_\ell)$, and the $\ell$-th power map $u\mapsto v=u^\ell$ induces $\mathrm{Tate}(q)/\mu_\ell\cong\mathrm{Tate}(q^\ell)$ with $\varphi^*(dv/v)=\ell\,du/u$ (see \cite[Ch.~8--10]{KatzMazur}, \cite[Ch.~V]{SilvermanATAEC}). Hence $\omega'=\ell^{-1}dv/v$ and
\[
\Delta^{[\ell]}\bigl(\mathrm{Tate}(q),\mu_\ell,\omega_{\mathrm{can}}\bigr)
=\ell^{-12}\,\Delta\bigl(\mathrm{Tate}(q^\ell),\ell^{-1}\tfrac{dv}{v}\bigr)
=\ell^{-12}\cdot\ell^{12}\,\Delta(q^\ell)=\Delta(q^\ell).
\]
The reduction of the classical $V_\ell(\Delta)$ is a Katz form of the same weight and level with the same $q$-expansion, and the $q$-expansion map is injective on forms of a given weight over the geometrically connected curve $X_0(\ell)_{\Fpb}$ \cite[\S1.6]{Katz350}; hence the two agree. Dividing by $\Delta$ and using that the completion at the cusp embeds the function field of the irreducible curve $X_0(\ell)_{\Fpb}$ into $\Fpb((q))$ gives the function-field statement; for $\ell=2$, Lemma~\ref{lem:identities} yields $R_2=\mu^2\eps/(64\mu\eps^2)=x/64$.
\end{proof}

In particular, for a non-cuspidal point $\wt P\in X_0(\ell)(\Fpb)$ represented by a pair $(E,C)$ and any nonzero $\omega$,
\begin{equation}\label{eq:Rvalue}
R_\ell(\wt P)\ =\ \frac{\ell^{-12}\,\Delta(E',\omega')}{\Delta(E,\omega)},\qquad \varphi^*\omega'=\omega,
\end{equation}
a quantity independent of $\omega$ and of the chosen representative.

\subsection{Deuring models}
The next step is to choose models on which Frobenius acts as the scalar $[-p]$. This normalization makes the subgroup and quotient in a supersingular isogeny rational over $\Fpp$.

\begin{lemma}\label{lem:unique-subgroup}
Let $E$ be a supersingular elliptic curve over an algebraically closed field $k$ of characteristic $p$. Then $E$ has a unique subgroup scheme of order $p$, namely $\ker F$, where $F\colon E\to E^{(p)}$ is the relative Frobenius and $E^{(p)}$ denotes the base change of $E$ along $x\mapsto x^p$; and a unique connected subgroup scheme of order $p^2$, namely $\ker(F^{(p)}\circ F)$, where $F^{(p)}\colon E^{(p)}\to E^{(p^2)}$ is the relative Frobenius of $E^{(p)}$.
\end{lemma}

\begin{proof}
Let $G\subset E$ be a subgroup scheme of order $p$. Since $E$ is supersingular, $E[p](k)=0$ \cite[Thm.~V.3.1]{SilvermanAEC}, so $G(k)=0$ and $G$ is connected (a nontrivial \'etale part would contribute $k$-points). A nontrivial connected finite group scheme has nonzero tangent space at the identity: its coordinate ring is a local $k$-algebra with $\mathfrak m\ne0$, so $\mathfrak m/\mathfrak m^2\ne0$ by Nakayama. The relative Frobenius $F_G$ of $G$ kills this tangent space, so $\ker F_G$ is a nontrivial subgroup scheme of $G$ of order dividing $p$; hence $\ker F_G=G$, i.e.\ $F_G=0$. By functoriality of Frobenius, $F_E$ restricted to $G$ equals $F_G=0$, so $G\subseteq\ker F_E$; both have order $p$, so $G=\ker F_E$.

Now let $G$ be connected of order $p^2$. As above, $\ker F_G=G\cap\ker F_E$ is nontrivial, hence equals $\ker F_E$ (the unique order-$p$ subgroup), so $\ker F_E\subset G$. Then $G/\ker F_E$ is a connected order-$p$ subgroup scheme of $E/\ker F_E\cong E^{(p)}$, which is again supersingular, so by the first part it equals $\ker F_{E^{(p)}}$; hence $G=\ker(F^{(p)}\circ F)$.
\end{proof}

\begin{lemma}[Deuring models]\label{lem:deuring}
Let $p\ge5$, let $E_0$ be a supersingular elliptic curve over $\Fpb$, and let $C_0\subset E_0$ be a subgroup of order $\ell\ne p$. Then there is a pair $(E,C)$, isomorphic to $(E_0,C_0)$ over $\Fpb$, such that:
\begin{enumerate}
\item[(a)] We have that $E$ is defined over $\Fpp$ and its $p^2$-power Frobenius endomorphism is $\pi_E=[-p]$;
\item[(b)] We have that $C$ is defined over $\Fpp$, hence so are $E':=E/C$ and the quotient isogeny $\varphi\colon E\to E'$;
\item[(c)] We have that $\pi_{E'}=[-p]$, and $E'$ is supersingular.
\end{enumerate}
\end{lemma}

\begin{proof}
(a) One has $j(E_0)\in\Fpp$ \cite[Thm.~V.3.1]{SilvermanAEC}; choose a model $E_1/\Fpp$ with this $j$-invariant, necessarily supersingular. Both $\pi_1:=\pi_{E_1}$ and $[p]$ are purely inseparable of degree $p^2$: for $\pi_1$ this is standard, and $[p]$ has degree $p^2$ with $E_1[p](\Fpb)=0$, so $E_1[p]=\ker[p]$ is connected and $[p]$ has trivial separable part. (By Lemma~\ref{lem:unique-subgroup}, the two kernels moreover coincide, both being the unique connected subgroup scheme of order $p^2$.) By \cite[Cor.~II.2.12]{SilvermanAEC}, each of $\pi_1$ and $[p]$ factors as a separable map following the $p^2$-power relative Frobenius $F^2\colon E_1\to E_1^{(p^2)}$: writing $\pi_1=\lam_1\circ F^2$ and $[p]=\lam_2\circ F^2$, the maps $\lam_1,\lam_2$ are separable isogenies of degree $\,p^2/p^2=1$, hence isomorphisms. Setting $u:=\lam_1\circ\lam_2^{-1}\in\Aut(E_{1,\Fpb})$ gives
\[
\pi_1=u\circ[p];
\]
moreover $u=\pi_1\circ[p]^{-1}$ in $\End(E_{1,\Fpb})\otimes\Q$ is fixed by $\Gal(\Fpb/\Fpp)$, since $\pi_1$ and $[p]$ are.

Now twist. Let $n:=\#\Aut(E_{1,\Fpb})\in\{2,4,6\}$. Since $n\mid 24\mid p^2-1$, we have $\mu_n\subset\Fpp$ and all automorphisms of $E_1$ are defined over $\Fpp$ (in a Weierstrass model they are $(x,y)\mapsto(\zeta^2x,\zeta^3y)$ with $\zeta\in\mu_n$), so $\Aut(E_{1,\Fpb})\cong\mu_n$ carries the trivial Galois action. The twists of $E_1$ over $\Fpp$ are classified by
\[
H^1\bigl(\Gal(\Fpb/\Fpp),\Aut(E_{1,\Fpb})\bigr)=\operatorname{Hom}_{\mathrm{cont}}(\widehat{\Z},\mu_n)\cong\mu_n,
\]
the class of a twist being determined by the value of its cocycle at $\sigma:=\Frob_{p^2}$ \cite[\S X.2, \S X.5]{SilvermanAEC}. For $\zeta\in\mu_n=\Aut(E_1)$ let $E_\zeta/\Fpp$ be the twist whose class sends $\sigma\mapsto\zeta$, with an $\Fpb$-isomorphism $\psi\colon E_1\to E_\zeta$ satisfying $\psi^\sigma=\psi\circ\zeta$ (here $\psi^\sigma$ is the conjugate isomorphism, obtained by applying $\sigma$ to the coefficients of $\psi$). For $P\in E_\zeta(\Fpb)$ with $Q:=\psi^{-1}(P)$, the arithmetic Frobenius acts by
\[
\sigma(P)=\sigma\bigl(\psi(Q)\bigr)=\psi^\sigma\bigl(\sigma Q\bigr)=\psi\bigl(\zeta(\pi_1 Q)\bigr).
\]
Since the Frobenius endomorphism of $E_\zeta$ induces the arithmetic Frobenius on $\Fpb$-points, and endomorphisms are determined by their action on points, $\pi_{E_\zeta}=\psi\circ\zeta\circ\pi_1\circ\psi^{-1}=(\psi\,\zeta u\,\psi^{-1})\circ[p]$. Choosing $\zeta:=-u^{-1}$ gives $\pi_{E_\zeta}=[-1]\circ[p]=[-p]$. Set $E:=E_\zeta$.

(b) Since $j(E)=j(E_0)$ there is an $\Fpb$-isomorphism $E_0\cong E$ \cite[Prop.~III.1.4]{SilvermanAEC}; transport $C_0$ to $C\subset E[\ell]$. The subgroup scheme $C$ is \'etale ($\ell\ne p$), so it is determined by $C(\Fpb)$, on which the Galois group acts through $\pi_E=[-p]$, i.e.\ by the scalar $-p$; a scalar preserves every subgroup, so $C$ is Galois-stable, hence defined over $\Fpp$, and then $E'=E/C$ and $\varphi$ are defined over $\Fpp$ \cite{Velu}.

(c) Since $\varphi$ is defined over $\Fpp$, one has $\pi_{E'}\circ\varphi=\varphi\circ\pi_E=\varphi\circ[-p]=[-p]\circ\varphi$, and $\varphi$ is an epimorphism, so $\pi_{E'}=[-p]$. Being isogenous to $E$, the curve $E'$ is supersingular.
\end{proof}

Once a Deuring model has been chosen, the value of $R_\ell$ is reduced to a ratio of the same discriminant invariant on the source and target curves.

\begin{lemma}\label{lem:kappa-reduction}
For an elliptic curve $E/\Fpp$ and a nonzero $\Fpp$-rational invariant differential $\omega$, the quantity
\[
\kappa(E)\ :=\ \Delta(E,\omega)^{(p^2-1)/12}\in\Fpp^\times
\]
is independent of the choice of $\omega$. If $\wt P\in X_0(\ell)(\Fpb)$ is represented by a pair $(E,C)$ as in Lemma~\textup{\ref{lem:deuring}}, then
\[
R_\ell(\wt P)^{\,m}\ =\ \frac{\kappa(E')}{\kappa(E)},\qquad m=\frac{p^2-1}{12}.
\]
\end{lemma}

\begin{proof}
Any two rational nonzero differentials differ by a scalar $\lam\in\Fpp^\times$, and
\[
\Delta(E,\lam\omega)^{m}
=\lam^{-12m}\Delta(E,\omega)^m
=\lam^{-(p^2-1)}\Delta(E,\omega)^m
=\Delta(E,\omega)^m.
\]
For the second claim, evaluate \eqref{eq:Rvalue} with a rational $\omega$. Choose a rational nonzero $\omega_0'$ on $E'$; since $\varphi$ is separable and rational, $\varphi^*\omega_0'=a\,\omega$ for some $a\in\Fpp^\times$, so the normalized differential is $\omega'=a^{-1}\omega_0'$ and $\Delta(E',\omega')=a^{12}\,\Delta(E',\omega_0')$. Therefore
\[
R_\ell(\wt P)^m=\ell^{-12m}\,a^{12m}\,\frac{\Delta(E',\omega_0')^m}{\Delta(E,\omega)^m}
=\ell^{-(p^2-1)}\,a^{\,p^2-1}\,\frac{\kappa(E')}{\kappa(E)}=\frac{\kappa(E')}{\kappa(E)},
\]
since $\ell\in\F_p^\times\subset\Fpp^\times$ and $a\in\Fpp^\times$.
\end{proof}

\subsection{Eta powers}
We prove the discriminant rationality theorem using roots of $\Delta$ that exist naturally on full level $3$ and full level $4$. The two levels together determine the required twelfth-root information.

\begin{proposition}\label{prop:eta}
For $N\in\{3,4\}$ set $k_N:=12/N$ (so $k_3=4$, $k_4=3$) and
\[
h_N(\tau):=\eta(\tau)^{2k_N}=q^{1/N}\prod_{n\ge1}(1-q^n)^{24/N},\qquad q^{1/N}:=e^{2\pi i\tau/N}.
\]
Then the following are true. 
\begin{enumerate}
\item[(a)] We have that $h_N\in M_{k_N}(\Gamma(N))$. Namely, for all $\gamma=\left(\begin{smallmatrix}a&b\\c&d\end{smallmatrix}\right)\in\Gamma(N)$ one has $h_N(\gamma\tau)=(c\tau+d)^{k_N}h_N(\tau)$, and $h_N$ is holomorphic at every cusp; moreover $h_N^{\,N}=\Delta$.
\item[(b)] Let $R_N:=\Z[\zeta_N,1/N]$ and let $\mathcal M_N$ denote the compactified moduli of elliptic curves equipped with a basis of the $N$-torsion whose Weil pairing \cite[\S III.8]{SilvermanAEC} equals $\zeta_N$; it is a smooth proper geometrically connected curve over $R_N$ \cite[Ch.~3--10]{KatzMazur}. Then $h_N$ arises from a unique algebraic modular form of weight $k_N$ over $R_N$ on $\mathcal M_N$, with $q$-expansion $q^{1/N}\prod_{n\ge1}(1-q^n)^{24/N}\in\Z[[q^{1/N}]]$ at the standard cusp, and $h_N^{\,N}=\Delta$ as algebraic forms. In particular $h_N(E,\phi,\omega)\ne0$ for every triple over a field.
\end{enumerate}
\end{proposition}

\begin{proof}
(a) Since $\eta^{24}=\Delta$ we have $h_N^N=\Delta$, and for $\gamma\in\Gamma(N)$ the quantity
\[
\chi(\gamma):=\frac{h_N(\gamma\tau)}{(c\tau+d)^{k_N}\,h_N(\tau)}
\]
satisfies $\chi(\gamma)^N=\Delta(\gamma\tau)\bigl((c\tau+d)^{12}\Delta(\tau)\bigr)^{-1}=1$; being a holomorphic function of $\tau$ with values in the finite set $\mu_N$, it is constant, and the cocycle property of $(c\tau+d)$ makes $\chi\colon\Gamma(N)\to\mu_N$ a group homomorphism. We must show $\chi=1$.

For $N\in\{3,4\}$: $-I\notin\Gamma(N)$, and $\Gamma(N)$ has no elliptic elements \cite[Ch.~2--3]{DiamondShurman}, so $\Gamma(N)$ acts freely and properly discontinuously on the upper half-plane $\mathbb H$ and is isomorphic to the fundamental group of $Y(N)=\mathbb H/\Gamma(N)$, a sphere with $c_N$ punctures ($c_3=4$, $c_4=6$; the genus is $0$) \cite[Ch.~3]{DiamondShurman}. The fundamental group of a punctured sphere is generated by loops encircling the punctures, and such a loop, lifted to $\mathbb H$, corresponds to a deck transformation preserving a small horoball at the corresponding cusp, hence to a parabolic element. Therefore $\Gamma(N)$ is generated by its parabolic elements.

Let $P\in\Gamma(N)$ be parabolic, fixing the cusp $s=\sigma\infty$ with $\sigma\in\SL_2(\Z)$. Then $\sigma^{-1}P\sigma$ is a parabolic element of $\SL_2(\Z)$ fixing $\infty$, so $\sigma^{-1}P\sigma=\pm T^{k}$ where $T=\left(\begin{smallmatrix}1&1\\0&1\end{smallmatrix}\right)$. Reducing modulo $N$: since $P\equiv I\pmod N$, also $\sigma^{-1}P\sigma\equiv I\pmod N$, because $\sigma^{-1}P\sigma-I=\sigma^{-1}(P-I)\sigma$ has entries in $N\Z$. The sign $-$ is impossible, since $-T^{k}\equiv I\pmod N$ would give $-1\equiv1\pmod N$, false for $N\ge3$; and $T^{k}\equiv I\pmod N$ forces $N\mid k$. Thus $P=\sigma T^{Nk'}\sigma^{-1}$ for some $k'\in\Z$.

Finally we compute $\chi(P)$. Set $g:=h_N|_{k_N}\sigma$. Then $g^N=\Delta|_{12}\sigma=\Delta$, and the reference function $g_0(\tau):=q^{1/N}\prod_{n\ge1}(1-q^n)^{24/N}$ also satisfies $g_0^N=\Delta$; both are holomorphic and nonvanishing on the connected set $\mathbb H$, so $g/g_0$ is a holomorphic function with values in $\mu_N$, i.e.\ a constant $\zeta\in\mu_N$: $g=\zeta g_0$. Since $g_0(\tau+N)=e^{2\pi i}\,q^{1/N}\prod_{n\ge1}(1-q^n)^{24/N}=g_0(\tau)$, we get $g|_{k_N}T^{Nk'}=g$, hence
\[
h_N|_{k_N}P=\bigl((h_N|\sigma)\,|\,T^{Nk'}\bigr)\,|\,\sigma^{-1}=g\,|\,\sigma^{-1}=h_N ,
\]
so $\chi(P)=1$ for every parabolic $P$, and therefore $\chi\equiv1$. Holomorphy at the cusps follows from the expansions $g=\zeta g_0$, whose $q^{1/N}$-exponents are nonnegative.

(b) The curve $\mathcal M_N\otimes\C$ is connected, the classical space $M_{k_N}(\Gamma(N),\C)$ is identified with the global sections of $\omega^{\otimes k_N}$ on it, and the algebraic $q$-expansion of a section at the standard cusp equals the analytic expansion in $q^{1/N}$. By Katz's $q$-expansion principle \cite[\S1.6]{Katz350} (see also \cite[Ch.~1]{KatzMazur}), a section over $\C$ whose $q$-expansion at a cusp of the geometrically connected curve has coefficients in $R_N$ descends uniquely to $R_N$; this applies to $h_N$, whose coefficients lie in $\Z$. The identity $h_N^N=\Delta$ persists over $R_N$ because both sides have equal $q$-expansions and the $q$-expansion map is injective (loc.\ cit.). Nonvanishing: $h_N(E,\phi,\omega)^N=\Delta(E,\omega)\ne0$.
\end{proof}

\subsection{Proofs of Theorems \ref{thm:kappa} and \ref{thm:containment}}
We now combine the eta-power calculation with the Frobenius action on torsion. The result is first the invariant $\kappa(E)$, and then the containment statement for $R_\ell$.

\begin{proof}[Proof of Theorem~\textup{\ref{thm:kappa}}]
Fix $N\in\{3,4\}$ and a ring homomorphism $R_N\to\Fpb$. Since $N\mid 24\mid p^2-1$, the residue degree of $p$ in $\Q(\zeta_N)$ --- the multiplicative order of $p$ modulo $N$ --- divides $2$, so the image of $R_N$ lies in $\Fpp$; hence the reduction of the algebraic form $h_N$ of Proposition~\ref{prop:eta} is defined over $\Fpp$, and its coefficients are fixed by $\sigma:=\Frob_{p^2}$. Let $\bar\zeta\in\mu_N(\Fpp)$ be the image of $\zeta_N$.

The group scheme $E[N]$ is \'etale, since $p\nmid N$. Choose a basis $\phi=(P_1,P_2)$ of $E[N](\Fpb)$ with Weil pairing $e_N(P_1,P_2)=\bar\zeta$ (possible, since the pairing takes every primitive value as bases vary), so that $(E,\phi)$ is an $\Fpb$-point of $\mathcal M_N$, and set
\[
v\ :=\ h_N(E,\phi,\omega)\in\Fpb^{\times}.
\]
The Galois action on torsion points of the $\Fpp$-curve $E$ is induced by its Frobenius endomorphism, so $\sigma\circ\phi=\pi_E\circ\phi=[-p]\circ\phi$; note that $[-p]$ acts on $E[N]$ as the scalar $-p\bmod N$ and multiplies the Weil pairing by $\det[-p]=p^2\equiv1\pmod N$, so $(E,\sigma\phi)$ is again an $\Fpb$-point of $\mathcal M_N$. Since $E$, $\omega$, and the coefficients of $h_N$ are all fixed by $\sigma$, Galois equivariance of evaluation gives
\[
\sigma(v)\ =\ h_N\bigl(E,\ [-p\bmod N]\cdot\phi,\ \omega\bigr).
\]
Now observe that $-p\bmod 3\in\{1,2\}=\{\pm1\}$ and $-p\bmod 4\in\{1,3\}=\{\pm1\}$ for \emph{every} prime $p\geq 5$. If $-p\equiv1\pmod N$ then $\sigma(v)=v$. If $-p\equiv-1\pmod N$, apply the isomorphism $[-1]\in\Aut(E)$: since $[-1]\circ\phi=-\phi$ and $[-1]^*\omega=-\omega$, isomorphism invariance and the weight rule give
\[
h_N(E,-\phi,\omega)=h_N\bigl(E,\phi,[-1]^*\omega\bigr)=h_N(E,\phi,-\omega)=(-1)^{k_N}\,h_N(E,\phi,\omega),
\]
so $\sigma(v)=(-1)^{k_N}v$.

\smallskip
\noindent
\emph{Case $N=3$, $k_3=4$:} in both cases $\sigma(v)=v$, so $v\in\Fpp^\times$ and, since $v^3=\Delta(E,\omega)$,
\[
\Delta(E,\omega)^{(p^2-1)/3}=v^{\,p^2-1}=1 .
\]

\smallskip
\noindent
\emph{Case $N=4$, $k_4=3$:} here $v^4=\Delta(E,\omega)$. If $p\equiv3\pmod4$ then $-p\equiv1\pmod4$, so $\sigma(v)=v$ and $\Delta(E,\omega)^{(p^2-1)/4}=v^{p^2-1}=1$. If $p\equiv1\pmod4$ then $-p\equiv-1\pmod4$, so $v^{p^2}=\sigma(v)=-v$ and $\Delta(E,\omega)^{(p^2-1)/4}=v^{\,p^2-1}=-1$.

Writing $\kappa:=\kappa(E)=\Delta(E,\omega)^{(p^2-1)/12}$, we have shown
\[
\kappa^4=\Delta(E,\omega)^{(p^2-1)/3}=1,\qquad
\kappa^3=\Delta(E,\omega)^{(p^2-1)/4}=
\begin{cases}
+1,& p\equiv3\ (4),\\
-1,& p\equiv1\ (4),
\end{cases}
\]
whence $\kappa=\kappa^4\cdot\kappa^{-3}=(\kappa^3)^{-1}=\kappa^3=(-1)^{(p+1)/2}$.
\end{proof}

\begin{proof}[Proof of Theorem~\textup{\ref{thm:containment}}]
Represent $\wt P$ by a pair $(E,C)$ as in Lemma~\ref{lem:deuring}. Both $E$ and $E'=E/C$ are supersingular over $\Fpp$ with Frobenius $[-p]$, so by Lemma~\ref{lem:kappa-reduction} and Theorem~\ref{thm:kappa},
\[
R_\ell(\wt P)^m=\frac{\kappa(E')}{\kappa(E)}=\frac{(-1)^{(p+1)/2}}{(-1)^{(p+1)/2}}=1. \qedhere
\]
\end{proof}

\begin{corollary}\label{cor:xvalues}
For every supersingular point $\wt P$ of $X_0(2)_{\Fpb}$, the coordinate $x=64R_2$ of Lemma~\textup{\ref{lem:hauptmodul}} satisfies
\[
x(\wt P)^{\,m}=64^{\,m},\qquad\text{hence}\qquad x(\wt P)^{\,t'}=64^{\,t'}\quad (t'=rm).
\]
In the notation of Proposition~\textup{\ref{prop:Apfactor}}: $\xi_i^{\,t'}=64^{\,t'}$ for $i=1,\dots,s$, and if $\alpha_p=1$ then $(-1)^{t'}=64^{\,t'}$.
\end{corollary}

\begin{proof}
By Theorem~\ref{thm:containment} at $\ell=2$ and Lemma~\ref{lem:Vell}, $(x(\wt P)/64)^m=R_2(\wt P)^m=1$; raise to the $r$-th power for the second equality. The last sentence is the identification of the supersingular $x$-values in Proposition~\ref{prop:Apfactor}.
\end{proof}

\begin{corollary}\label{cor:classical}
Let $E/\Fpp$ be supersingular with $\pi_E=[-p]$ and let $\omega$ be rational. Then $\Delta(E,\omega)$ is both a square and a cube in $\Fpp^\times$. Consequently every supersingular $j$-invariant is a cube in $\Fpp$, and $j-1728$ is a square in $\Fpp$.
\end{corollary}

\begin{proof}
By Theorem~\ref{thm:kappa}, $\Delta(E,\omega)^{(p^2-1)/3}=\kappa^4=1$ and $\Delta(E,\omega)^{(p^2-1)/2}=(\kappa^3)^2=1$. For a rational Weierstrass model one has $c_4,c_6\in\Fpp$ with $j=c_4^3/\Delta$ and $j-1728=c_6^2/\Delta$ \cite[\S III.1]{SilvermanAEC}; the cases $j=0$ and $j=1728$ are trivial, and every supersingular $j$ arises from a Deuring model by Lemma~\ref{lem:deuring}.
\end{proof}

\begin{remark}\label{rem:sharp}
The exponent $12$ is sharp in two senses. First, for $p\equiv1\pmod4$ Theorem~\ref{thm:kappa} says $\kappa(E)=-1$: the discriminant of a Deuring model is a sixth power but \emph{not} a twelfth power in $\Fpp^\times$. Second, the supersingular values of $R_\ell$ need not be $24$th powers in $\Fpp^\times$. For example, at $p=11$ the supersingular $x$-values on $X_0(2)$ are $-1,-3,-4$ (Section~\ref{sec:examples}), and $64^{-1}\equiv5\pmod{11}$, so $R_2$ takes the values $6,7,2$ at the three supersingular points; each satisfies $R_2^{10}=1$ (here $m=10$), in accordance with Theorem~\ref{thm:containment}, but
\[
6^{(p^2-1)/24}=6^5=-1,\qquad 7^5=-1,\qquad 2^5=-1\qquad\text{in }\F_{11}\subset\F_{121}.
\]
Hence none of these values is a $24$th power in $\F_{121}^\times$, and no $\F_{p^2}$-rational $24$th root of $R_2$ exists at any supersingular point of $X_0(2)_{\overline{\F}_{11}}$. In particular, although the eta quotient $g_\ell=\eta(\ell z)/\eta(z)$ satisfies $g_\ell^{24}=R_\ell$, its supersingular values do \emph{not} in general lie in $\Fpp$ (at $p=11$ they generate $\F_{11^4}$): $g_\ell$ lives on a degree-$24$ Kummer covering of $X_0(\ell)$ whose residue fields exceed $\Fpp$. Any proposed proof of Theorem~\ref{thm:containment} that descends $\Delta^{1/24}$ or $\Delta^{1/12}$ to $\Fpp$ must therefore fail; the argument above succeeds precisely because it descends only $\Delta^{1/3}$ and $\Delta^{1/4}$, for which the Frobenius scalar $-p$ acts on the relevant torsion as $\pm1$.
\end{remark}

\section{Proof of Theorem \ref{thm:main}}\label{sec:proof}

We now assemble the arithmetic results into Behrens' four conditions. The level-$2$ factorization gives the sharp bound, while the containment theorem gives the divisibility at every level.

\begin{lemma}[Uniform multiplicity]\label{lem:LucasKummer}
Let $k$ be a field of characteristic $p$, let $b\in k^{\times}$, and let $t\in\Z_{\ge1}$ with $\nu_p(t)=n$, so $t=p^{n}t'$ with $p\nmid t'$. Then, over an algebraic closure $\bar k$,
\[
x^{t}-b\ =\ \bigl(x^{t'}-b_0\bigr)^{p^{n}}\ =\ \prod_{\rho\in W}(x-\rho)^{p^n},\qquad b_0:=b^{\,p^{-n}}\in\bar k^{\times},
\]
where $W:=\{\rho\in\bar k^{\times}:\rho^{t'}=b_0\}$ has exactly $t'$ elements. In particular $x^{t'}-b_0$ is separable and every root of $x^{t}-b$ has multiplicity exactly $p^{n}$.
\end{lemma}

\begin{proof}
Since $y\mapsto y^{p}$ is an automorphism of $\bar k$, there is a unique $b_0=b^{p^{-n}}\in\bar k^{\times}$ with $b_0^{p^n}=b$, and
\[
x^{t}-b=(x^{t'})^{p^{n}}-b_0^{\,p^{n}}=\bigl(x^{t'}-b_0\bigr)^{p^{n}}
\]
because raising to the $p^{n}$-th power is additive in characteristic $p$. The derivative of $x^{t'}-b_0$ is $t'x^{t'-1}\ne0$, whose only root is $x=0$; as $b_0\ne 0$, no root of $x^{t'}-b_0$ is $0$, so $x^{t'}-b_0$ is separable with exactly $t'$ roots, all nonzero.
\end{proof}

The uniform multiplicity lemma lets us separate the $p^n$-power contribution from the prime-to-$p$ exponent $t'$.

\begin{proposition}[Factorization of $L_2(\Delta^t)$]\label{prop:L2factor}
In $\Fpb[\delta,\eps]$, we have
\[
L_2(\Delta^t)\ =\ \mu^{t}\eps^{t}\bigl(\mu^{t}-64^{t}\eps^{t}\bigr)
\ =\ \mu^{t}\eps^{t}\prod_{\rho\in W}(\mu-\rho\,\eps)^{p^n},
\qquad W:=\{\rho\in\Fpb^{\times}:\rho^{\,t'}=64^{\,t'}\},
\]
with $\#W=t'$ and $0\notin W$. If $-1\in W$, the corresponding factor is $(\mu+\eps)^{p^n}=\delta^{2p^n}$.
\end{proposition}

\begin{proof}
Since $V_2$ is a ring homomorphism, Lemma~\ref{lem:identities} gives, modulo $p$,
\[
L_2(\Delta^t)=V_2(\Delta)^{t}-\Delta^{t}=\mu^{2t}\eps^{t}-64^{t}\mu^{t}\eps^{2t}
=\mu^{t}\eps^{t}\bigl(\mu^{t}-64^{t}\eps^{t}\bigr).
\]
Apply Lemma~\ref{lem:LucasKummer} with $b=64^{t}$; then $b_0=64^{t'}$ (indeed $(64^{t'})^{p^n}=64^{t}$), so
$x^{t}-64^{t}=\prod_{\rho\in W}(x-\rho)^{p^n}$. Substituting $x=\mu/\eps$ and clearing denominators --- i.e.\ multiplying by $\eps^{t}=\prod_{\rho\in W}\eps^{p^n}$, which is legitimate as $t=t'p^n=(\#W)p^n$ --- homogenizes this to $\mu^{t}-64^{t}\eps^{t}=\prod_{\rho\in W}(\mu-\rho\eps)^{p^n}$. Finally $\mu+\eps=\delta^{2}$.
\end{proof}

The factorization gives the exact level-$2$ divisibility by comparing its linear factors with those of $A_p$.

\begin{proposition}[Divisibility at level $2$]\label{prop:level2}
If $j\ge1$, then we have that
\[
A_p^{\,j}\ \bigm|\ L_2(\Delta^t)\ \text{ in }\ \F_p[\delta,\eps]
\qquad\Longleftrightarrow\qquad j\le p^{n}.
\]
\end{proposition}

\begin{proof}
Divisibility in $\F_p[\delta,\eps]$ is equivalent to divisibility in $\Fpb[\delta,\eps]$: if $L_2(\Delta^t)=A_p^{\,j}H$ with $H\in\Fpb[\delta,\eps]$, then applying any element of $\Gal(\Fpb/\F_p)$ and cancelling $A_p^{\,j}$ in the integral domain shows $H$ is Galois-invariant, hence lies in $\F_p[\delta,\eps]$.

Work in the unique factorization domain $\Fpb[\delta,\eps]$. The elements $\delta$, $\eps$, $\mu$, and $\mu-\rho\eps$ for $\rho\in\Fpb\smallsetminus\{0,-1\}$ are irreducible and pairwise non-associate (they are homogeneous of weight $2$ or $4$ and have distinct zero loci on $X_0(2)$ by Lemma~\ref{lem:hauptmodul}). By Corollary~\ref{cor:xvalues}, every $\xi_i$ lies in $W$; and if $\alpha_p=1$, then the point $x=-1$ is supersingular by Proposition~\ref{prop:Apfactor}, so $-1\in W$. From Proposition~\ref{prop:L2factor} we read off the valuations of $L:=L_2(\Delta^t)$ at the relevant primes:
\[
v_{\mu-\xi_i\eps}(L)=p^{n}\quad(1\le i\le s),\qquad
v_{\delta}(L)=\begin{cases}2p^{n},& -1\in W,\\ 0,&\text{else,}\end{cases}
\]
where for the first we use that $\xi_i\in W$ occurs exactly once in the product (the roots are distinct) and that $\xi_i\ne0$, so $\mu^{t}\eps^{t}$ contributes nothing; and for the second that modulo $\delta$ one has $\mu\equiv-\eps\ne0$ and $\mu-\rho\eps\equiv-(1+\rho)\eps\ne0$ for $\rho\ne-1$, while $\mu+\eps=\delta^2$. By Proposition~\ref{prop:Apfactor}, $A_p^{\,j}=c^{\,j}\delta^{\,\alpha_p j}\prod_{i=1}^{s}(\mu-\xi_i\eps)^{\,j}$.

If $j\le p^{n}$: each $v_{\mu-\xi_i\eps}(L)=p^{n}\ge j$; and if $\alpha_p=1$ then $-1\in W$ and $v_\delta(L)=2p^{n}\ge j$. Hence $A_p^{\,j}\mid L$.

If $j>p^{n}$: since $s\ge1$ (Proposition~\ref{prop:Apfactor}), the prime $\mu-\xi_1\eps$ divides $A_p^{\,j}$ to order $j$ but divides $L$ to order exactly $p^{n}<j$. Hence $A_p^{\,j}\nmid L$.
\end{proof}

\begin{remark}\label{rem:elliptic}
When $p\equiv3\pmod4$ the elliptic point $x=-1$ of $X_0(2)$ is supersingular and contributes the factor $\delta$ to $A_p$; but $L_2(\Delta^t)$ is divisible by $\delta^{2p^n}$ --- \emph{twice} the naive count --- because the coordinate ramifies there: $x+1=\delta^2/\eps$. Had $A_p$ consisted of the factor $\delta$ alone, the divisibility $A_p^{\,j}\mid L_2(\Delta^t)$ would persist up to $j=2p^n$ and clause (ii) of Theorem~\ref{thm:main} would fail at $\ell=2$. It is therefore essential that a non-elliptic supersingular factor $\mu-\xi_1\eps$ exists; the count $s=(p-1-2\alpha_p)/4\ge1$ of Proposition~\ref{prop:Apfactor} guarantees this for every $p\ge5$.
\end{remark}

For the positive direction of the conjecture, the same divisibility must hold at every prime level. This is where the supersingular containment theorem replaces the explicit level-$2$ calculation.

\begin{proposition}[Divisibility at every level]\label{prop:bridge}
For every prime $\ell\ne p$ and every $1\le j\le p^{n}$, we have
\[
A_p^{\,j}\ \bigm|\ L_\ell(\Delta^t)\qquad\text{in } M_*(\ell)_{\Z/p}.
\]
\end{proposition}

\begin{proof}
Since $V_\ell$ is a ring homomorphism and we are in characteristic $p$,
\begin{equation}\label{eq:freshman}
L_\ell(\Delta^t)=\bigl(V_\ell\Delta\bigr)^{t}-\Delta^{t}
=\Bigl(\bigl(V_\ell\Delta\bigr)^{t'}\Bigr)^{p^n}-\bigl(\Delta^{t'}\bigr)^{p^n}
=\Bigl(L_\ell\bigl(\Delta^{t'}\bigr)\Bigr)^{p^n}.
\end{equation}
It therefore suffices to prove $A_p\mid G$ for $G:=\overline{L_\ell(\Delta^{t'})}$, since then $A_p^{\,p^n}\mid G^{\,p^n}=\overline{L_\ell(\Delta^t)}$ and a fortiori $A_p^{\,j}\mid\overline{L_\ell(\Delta^t)}$ for $j\le p^n$.

First, $G$ vanishes at every supersingular point of $X_0(\ell)_{\Fpb}$: by Lemma~\ref{lem:Vell}, the value of $G$ on a triple $(E,C,\omega)$ over a field is
\[
\Delta^{[\ell]}(E,C,\omega)^{t'}-\Delta(E,\omega)^{t'}
=\Delta(E,\omega)^{t'}\bigl(R_\ell(\wt P)^{t'}-1\bigr),
\]
and at a supersingular $\wt P$ one has $R_\ell(\wt P)^{t'}=\bigl(R_\ell(\wt P)^{m}\bigr)^{r}=1$ by Theorem~\ref{thm:containment}.

Now divide. Choose $N\in\{3,4\}$ with $\gcd(N,\ell)=1$ ($N=3$ unless $\ell=3$, in which case $N=4$), and let $\overline Y$ be the compactified moduli, over $\Fpb$, of elliptic curves equipped with a $\Gamma_0(\ell)$-structure and a full level-$N$ structure of fixed Weil pairing. This is a smooth proper geometrically connected curve (the level-$N$ moduli problem is representable and rigid for $N\ge3$, and adding the $\Gamma_0(\ell)$-structure preserves this; see \cite[Ch.~3--8]{KatzMazur}), equipped with an action of a finite group $Q$ (a subgroup of $\GL_2(\Z/N)$, of order prime to $p$ since $p\ge5$) such that the $Q$-invariant sections of $\omega^{\otimes k}$ on $\overline Y$ are the weight-$k$ elements of $M_*(\ell)_{\Z/p}\otimes\Fpb$. (Here $\omega$ denotes, as usual, the line bundle on the moduli whose $k$-th tensor power has the weight-$k$ forms as its sections \cite[\S1]{Katz350}; a rule on triples in the sense of \S\ref{subsec:hasse} is the same thing as such a section, by evaluation on the universal triple.) By Igusa's theorem, the pullback of $A_p$ to $\overline Y$ has \emph{simple} zeros, located exactly at the supersingular points \cite{Igusa}, \cite[Ch.~12]{KatzMazur}, and by Lemma~\ref{lem:hasse}(b) it is nonvanishing at the cusps. The pullback of $G$ vanishes at every supersingular point of $\overline Y$, as shown above. Therefore, we have that
\[
H:=\frac{G}{A_p}
\]
is a section of $\omega^{\otimes(12t'-(p-1))}$ on all of $\overline Y$: its only possible poles, at the simple zeros of $A_p$, are cancelled by the vanishing of $G$. Moreover $H$ is $Q$-invariant, being a ratio of $Q$-invariant sections. Hence $H\in M_{12t'-(p-1)}(\ell)_{\Z/p}\otimes\Fpb$ and $G=A_p\cdot H$. Finally $H=G/A_p$ is fixed by $\Gal(\Fpb/\F_p)$, since $G$ and $A_p$ are and quotients in a domain are unique; so $H$ is defined over $\F_p$ and the divisibility holds in $M_*(\ell)_{\Z/p}$.
\end{proof}

\begin{remark}\label{rem:rings}
Two remarks on the coefficient ring in (C4). Both rest on the following base-change facts for the fine moduli $\overline Y$ of the proof of Proposition~\ref{prop:bridge} (for $\ell=2$ take $N=3$). Let $\mathcal O$ be the localization of $R_N[1/\ell]$ at a prime above $p$ --- a discrete valuation ring of mixed characteristic $(0,p)$ whose fraction field embeds in $\C$ --- over which $\overline Y$ is a smooth proper curve with geometrically connected fibers, carrying the line bundle $\omega$ and the action of $Q$. By the Kodaira--Spencer isomorphism $\omega^{\otimes2}\cong\Omega^1_{\overline Y}(\mathrm{cusps})$ \cite[\S A1.3]{Katz350}, for every $k\ge2$ the fiberwise degree of $\omega^{\otimes k}$ equals $(2g-2)+c+(k-2)\deg\omega>2g-2$, where $g$ is the common genus of the fibers and $c\ge1$ the number of cusps; hence $H^1(\omega^{\otimes k})$ vanishes on both fibers, and by cohomology and base change $H^0(\overline Y,\omega^{\otimes k})$ is a free $\mathcal O$-module whose formation commutes with arbitrary base change \cite[\S1.7]{Katz350}. Since $p\nmid\#Q$, the projector $e_Q:=\tfrac1{\#Q}\sum_{g\in Q}g$ is defined over $\mathcal O$ and its image, the $Q$-invariants,  is a free direct summand, again compatible with base change. Two consequences follow for $k\ge2$. First, the dimension count
\begin{equation}\label{eq:dims}
\dim_{\Fpb}\bigl(M_k(\ell)_{\Z/p}\otimes\Fpb\bigr)\;=\;\dim_{\Fpb}H^0(\overline Y_{\Fpb},\omega^{\otimes k})^{Q}\;=\;\dim_{\C}H^0(\overline Y_{\C},\omega^{\otimes k})^{Q}\;=\;\dim_{\C}M_k(\Gamma_0(\ell),\C),
\end{equation}
the last equality because $Q$-invariant classical forms of level $\Gamma(N)\cap\Gamma_0(\ell)$ (on the component of fixed Weil pairing) are exactly the forms on $\Gamma_0(\ell)$. Second, every geometric form of weight $k$ lifts to characteristic zero: the reduction map $H^0(\overline Y,\omega^{\otimes k})\to H^0(\overline Y_{\Fpb},\omega^{\otimes k})$ is surjective by freeness and base change, and averaging a lift over $Q$ preserves the reduction while landing in the $Q$-invariants.

With these facts in hand: first, at $\ell=2$ the geometric ring coincides with Larson's ring \eqref{eq:ring} in every weight that occurs in this paper (the even weights $12t-j(p-1)\ge2$ and, trivially, weight $0$). Indeed, the reduction of a classical form is a geometric form, giving an inclusion $\F_p[\delta,\eps]_k\subseteq M_k(2)_{\Z/p}$, and by \eqref{eq:dims} the two sides have the common dimension
\[
\dim_{\C}M_k(\Gamma_0(2),\C)\;=\;\Bigl\lfloor\frac k4\Bigr\rfloor+1\;=\;\#\{(a,b)\in\Z_{\ge0}^2:2a+4b=k\}\;=\;\dim\F_p[\delta,\eps]_k
\]
for even $k\ge2$, so the inclusion is an equality. Hence Proposition~\ref{prop:level2} settles condition (C4) at level $2$ in either reading, and this is the only level used in clause (ii). Second, the divisibility of Proposition~\ref{prop:bridge} also holds in the classical reading, i.e.\ there is $h\in M_{12t-j(p-1)}(\Gamma_0(\ell))$ with $p$-integral coefficients and $L_\ell(\Delta^t)\equiv E_{p-1}^{\,j}\,h\pmod p$: the geometric form $A_p^{\,p^n-j}H^{\,p^n}$ has weight $12t-j(p-1)\ge2$, so by the lifting statement above it is the reduction of a $Q$-invariant form over $\mathcal O$, which descends to a classical form $h$ on $\Gamma_0(\ell)$ in characteristic zero with the required reduction.
\end{remark}

It remains only to check the three Behrens conditions that do not involve the operators $L_\ell$.

\begin{lemma}\label{lem:C123}
For every admissible pair $(i,j)$, the form $\Delta^t$ satisfies \textup{(C1)}, \textup{(C2)} and \textup{(C3)}.
\end{lemma}

\begin{proof}
(C1): the reduction of $\Delta^t$ has leading coefficient $1$, so it is nonzero. (C2): by Lemma~\ref{lem:maxord}, a nonzero reduction of a form of weight $12t$ has $q$-order at most $t$, and $\overline{\Delta^t}$ attains this maximum, $\ord_q(\overline{\Delta^t})=t$; hence, since $j\ge1$,
\[
12\ord_q(\overline{\Delta^t})=12t>12t-(p-1)j,
\]
which is the first alternative in (C2). (C3) is Lemma~\ref{lem:filtration}.
\end{proof}

\begin{proof}[Proof of Theorem~\textup{\ref{thm:main}}]
(i) Let $j$ be admissible with $1\le j\le p^{n}$. Conditions (C1)--(C3) hold for $\Delta^t$ by Lemma~\ref{lem:C123}, and (C4) holds at every prime $\ell\ne p$ by Proposition~\ref{prop:bridge}. Hence $\Delta^t$ satisfies all four of Behrens' conditions and may serve as $f_{i/j}$.

(ii) By Larson's enumeration \cite[\S2, Lem.~2.1]{Larson}, admissible indices $j>p^{n}$ exist precisely when $r>1$ and $n\ge2$: for $r=1$ the range ends at $p^n$, for $r>1,n=0$ the only admissible value is $j=1=p^0$, and for $r>1,n\ge1$ the range extends to $b_n=p^{n}+p^{n-1}-1$, which exceeds $p^n$ exactly when $n\ge2$ (moreover all excluded indices $p,2p,\dots,b_{n-2}p$ are $<p^n$, so every integer in $(p^n,b_n]$ is admissible). For such $j$, Proposition~\ref{prop:level2} shows $A_p^{\,j}\nmid L_2(\Delta^t)$ in $M_*(2)_{\Z/p}=\F_p[\delta,\eps]$. Hence $\Delta^t$ fails condition (C4) at $\ell=2$, so it cannot be a Behrens form, and no representative of $f_{i/j}$ equals $\Delta^t$.
\end{proof}

\section{Examples}\label{sec:examples}

These examples show how the abstract containment statement recovers Larson's explicit computations. We illustrate the mechanism at Larson's primes.

\subsection{The case of $p=11$}
This is the first case in which the elliptic point contributes a doubled factor. The non-elliptic supersingular factors still give the sharp bound.

Here $m=(11^2-1)/12=10$ and $64\equiv 9\pmod{11}$, with $9^5=1$ in $\F_{11}$, so $64^{t'}=64^{10r}=1$. Proposition~\ref{prop:Apfactor} gives $\alpha_p=1$ ($11\equiv3\bmod4$) and $s=2$; explicitly, up to a unit,
\[
A_{11}=\delta\,(\mu+3\eps)(\mu+4\eps),
\]
recovering \cite[\S5]{Larson}, so the supersingular $x$-values are $-1$ (the elliptic point), $-3$, and $-4$. Corollary~\ref{cor:xvalues} predicts $\xi^{10r}=1$ for each, which one checks by hand: $(-1)^{10r}=1$, while $(-3)^{10}=9^{5}=1$ and $(-4)^{10}=16^{5}=5^{5}=1$ in $\F_{11}$. For $t=10i$, $t'=10r$, Proposition~\ref{prop:L2factor} gives
\[
L_2(\Delta^{t})=\mu^{t}\eps^{t}\prod_{\rho^{10r}=1}(\mu-\rho\eps)^{11^{n}},
\]
and the roots $\rho\in\mu_{10r}(\Fpb)$ include all of $\F_{11}^\times=\mu_{10}$, in particular $-1,-3,-4$; the factor at $\rho=-1$ is $(\mu+\eps)^{11^n}=\delta^{2\cdot 11^{n}}$. Hence $A_{11}^{\,j}\mid L_2(\Delta^{t})$ exactly for $j\le 11^{n}$, the binding constraint coming from $(\mu+3\eps)$ or $(\mu+4\eps)$; the elliptic factor $\delta$ alone would permit $j\le 2\cdot 11^{n}$, illustrating Remark~\ref{rem:elliptic}. This matches Larson's computations \cite[\S5]{Larson}.

\subsection{The case of $p=13$}
Here the supersingular values include a conjugate pair over $\F_{13^2}$. The example illustrates why the argument is naturally formulated over $\Fpp$.

Here $m=14$ and $64\equiv-1\pmod{13}$, so again $64^{t'}=64^{14r}=1$. Now $\alpha_p=0$ ($13\equiv1\bmod4$) and $s=3$; up to a unit,
\[
A_{13}=(\mu+12\eps)\,(\mu^{2}+5\mu\eps+\eps^{2}),
\]
recovering \cite[\S5]{Larson}: the supersingular $x$-values are $1\ (=-12)$ together with the two roots $\alpha,\alpha^{13}$ of $x^{2}+5x+1$, which lie in $\F_{169}\smallsetminus\F_{13}$ (the ``inert'' case, illustrating that supersingular values live in $\Fpp$ rather than $\F_p$). Containment is visible by hand: $1^{14r}=1$; and since $\alpha\cdot\alpha^{13}=\mathrm{N}_{\F_{169}/\F_{13}}(\alpha)=1$ (the constant term of $x^2+5x+1$), we get $\alpha^{14}=\alpha^{13}\cdot\alpha=1$, so $\alpha^{14r}=1$ as well. Thus all three supersingular values lie in $\mu_{14}\subset W$, and $A_{13}^{\,j}\mid L_2(\Delta^{t})$ exactly for $j\le 13^{n}$, in agreement with \cite[\S5]{Larson}.

\end{document}